\documentclass[12pt]{amsart}

\usepackage{mathrsfs}
\usepackage{eufrak} 
\usepackage{amscd}
\usepackage{amssymb}
\usepackage{a4wide}
\usepackage{amstext}
\usepackage{amsthm}
\usepackage{color}
\usepackage[T2A]{fontenc}

\newcommand{\Z}{\mathbb{Z}}

\renewcommand{\phi}{\varphi}

\newcommand{\N}{\mathbb{N}}

\newcommand{\Lin}[1]{\overline{Lin} \{ #1\}}
\DeclareMathOperator*\lowlim{\underline{lim}}

\theoremstyle{plain}
\newtheorem{Theorem}{Theorem}

\newtheorem{Lemma}[Theorem]{Lemma}

\newtheorem{Corollary}[Theorem]{Corollary}

\theoremstyle{definition}
\newtheorem{Example}[Theorem]{Example}
\newtheorem{Definition}[Theorem]{Definition}
\newtheorem{Remark}[Theorem]{Remark}

\begin{document}
\sloppy

\title[]{Hereditarily and nonhereditarily complete systems of vectors in a Hilbert space }

\author{Mikhail Prokofyev}
\address{St. Petersburg State University}
\email{mikhail.prokofyev@yandex.ru}

%\date{}


\subjclass[]{}
\thanks{The work is supported by the Russian Science Foundation grant No. 23-11-00153.}

\begin{abstract}
In this paper, we study the property of hereditary completeness of vector systems $\{x_k\}_{k=1}^\infty$ in a Hilbert space. A criterion of hereditary completeness is obtained in terms of projectors on closed linear spans of systems of the form $\{x_k\}_{k \in N}$, $N \subset \N$. Developed technique has been used to prove that mixed systems of a hereditarily complete system are also hereditarily complete. In conclusion, the problem of possible defects in a nonhereditarily complete system is considered. 
\end{abstract}

\maketitle

%\maketitle

\section{Introduction}
\label{section1}
This paper is devoted to complete minimal systems in a separable Hilbert space, as well as their special classes. Let us introduce some notation and definitions.

Everywhere in what follows, $H$ will denote a separable Hilbert space. If $A \subset X$ is a subset of some set $X$, then $A^c$ will denote the complement of $A$ in $X$ ($X$ should be clear from the context wherever such a notation is used). If $A,B \subset X$ and $A \cap B=\emptyset$, then the union of $A$ and $B$ will be denoted as $A \sqcup B$. 
 \begin{Definition}
   Let $\{x_k\}_{k=1}^\infty$ be a system of vectors in $H$. It is called \textbf{minimal} if for any $k$ the vector $x_k $ does not lie in  $\Lin{x_j}_{j \in \N \setminus \{k\}}$.
 \end{Definition}
  It is well known (and follows immediately from the Hahn--Banach theorem) that minimality is equivalent to the following property: the system $\{x_k\}_{k=1}^\infty$ has a biorthogonal system $\{x_k^*\}_{k=1}^\infty$, that is, a vector system with the property $(x_k,x_j^*)=\delta_{k,j}$, where $\delta_{k,j}$ is the Kronecker delta symbol. Everywhere in the sequel, if $\{x_k\}_{k=1}^\infty$ is a system of vectors, then $\{x_k^*\}_{k=1}^\infty$ is a system biorthogonal to it. If the vector system is complete, then the biorthogonal one is unique, but in the general case there is no uniqueness.

 In connection with the spectral theory of linear operators, Marcus \cite{M} introduced the following concept of strong completeness which later became known as hereditary completeness \cite{DN}. In some sources, it is also called the strong Markushevich basis or the strong $M$--basis.
 \begin{Definition}
  Let $\{x_k\}_{k=1}^\infty$ be the complete system of vectors in $H$. It is called \textbf{hereditarily complete} if it is minimal and for each $x \in H$ the inclusion $x \in \Lin{(x,x_k^*)x_k}_{k = 1}^\infty $ holds. 
 \end{Definition}
  
 Hereditary completeness generalizes different concepts of a basis in a Hilbert space. Indeed, let any vector $x\in H$ decompose into a series $x=\sum\limits_{k=1}^\infty c_k x_k$, then multiplying by $x_k^*$, we get $(x,x_k^*)=c_k$ and $x \in\Lin{(x,x_k^*)x_k}_{k = 1}^\infty $.

It is not difficult to show that the definition of hereditary completeness is equivalent to the following property.
 \begin{Lemma}
  \label{hered_comp_crit}
The system $\{x_k\}_{k=1}^\infty$ is hereditarily complete if and only if for any partition into disjoint sets $\mathbb{N}=N_1\sqcup N_2$ the system $\{x_k\}_{k\in N_1} \cup\{x_k^*\}_{k\in N_2}$ is complete in $H$. 
\end{Lemma} 
In connection with this criterion, we introduce two definitions.

\begin{Definition}
    Systems of the form $\{x_k\}_{k\in N} \cup \{x_k^*\}_{k\in N^c}$ will be called \textbf{mixed systems} corresponding to the minimal system $\{x_k\}_{k=1}^\infty$ and the set $N \subset \N$. 
\end{Definition}

\begin{Definition} \label{def}
Let $\{x_k\}_{k=1}^\infty$ be a vector system that has a biorthogonal one, $N\subset\mathbb{N}$. We will say that the system has a \textbf{defect} $k$ with respect to $N$ ($k$ can be infinite) if the dimension of the orthogonal complement of $\overline{Lin}(\{x_k\}_{k\in N}\cup\{x_k^*\}_{k \in N^c })$ is equal to $k$.
\end{Definition}
Obviously, it follows from the property of hereditary completeness that the biorthogonal system is complete. In general, a system biorthogonal to a complete system does not have to be complete itself. The simplest example of such a situation is given by the system $\{e_1+e_k\}_{k\geq 2}$, where $\{e_k\}_{k\geq1}$ is an orthonormal basis. Moreover, Young \cite{RY} showed that an arbitrary closed infinite-dimensional subspace can be a closed linear span of a biorthogonal system (see Example \ref{e_1+e_k} and Lemma \ref{example_inf} below). 

Constructing a complete minimal system $\{x_k\}$ with a complete biorthogonal system $\{x_k^*\}$, but without the property of hereditary completeness (i.e., such that some nontrivial mixed system is incomplete), presents certain difficulties. There are several different construction methods:
\begin{enumerate}
    \item Markus' original example \cite{M};
    \item General constructions of nonhereditarily complete systems \cite{DN,DNS};
    \item Diagonal systems (i.e., such that $x_n\in Lin\{e_{n-L},\ldots,e_{n+L}\}$), which were studied in \cite{AS, KLP, P2020, P2022};
    \item Nonhereditarily complete systems of exponentials on the segment \cite{BBB1} and reproducing kernels in the spaces of entire functions \cite{BBB2,BBB3}.
\end{enumerate}
Most of the works cited above deal with the construction of a single incomplete mixed system. Moreover, it often follows from the construction that a possible defect in a mixed system does not exceed 1 (or takes all values from 0 to a certain number $n$, as in \cite{BBB2}). A natural question arises: is it possible to construct a nonhereditarily complete system that realizes exactly a given (finite or infinite) set of defects? One of the goals of this paper is to give a partial positive answer to this question (see  Theorem \ref{def_custom}).

We will also discuss a number of equivalent reformulations of the hereditary completeness property in terms of subspaces generated by vectors of the system, and in terms of convergence of some families of projectors (Corollary \ref{hered_comp_proj} and Lemma \ref{hered_comp_top}). As a consequence, we show that if a system is hereditarily complete, then all of its mixed systems are hereditarily complete (Theorem \ref{mix_hered_comp}).
\section{The simplest properties and examples of minimal systems}

For completeness, we present a proof of the Lemma \ref{hered_comp_crit}.
\begin{proof}
Let us first assume that the system $\{x_k\}_{k=1}^\infty$ is hereditarily complete. Let the vector $x$ be orthogonal to the system $\{x_k\}_{k\in N_1}\cup\{x_k^*\}_{k\in N_2}$. Since $x$ belongs to the space $\overline{Lin}\{ (x,x_k^*)x_k\}_{k=1}^\infty$, it also lies in $\Lin{x_k}_{k\in N_1}$. However, $x$ is orthogonal to the last subspace, and therefore is equal to zero.

To prove the converse, consider an arbitrary element $x$ in $H$. Let 
$$
N_2=\{ k \in \mathbb{N} \ | \ (x,x_k^*)=0\}
$$
and $N_1=N_2^c$. Since $(x_k,x_j^*)=0$ for different $k$ and $j$, we conclude from the completeness of the system $\{x_k\}_{k\in N_1} \cup \{x_k^*\}_{k\in N_2}$ that $\Lin{x_k}_{k \in N_1}$ is the orthogonal complement of $\Lin{x_k^*}_{k \in N_2}$. But then  $x$ lies in $\Lin{x_k}_{k \in N_1}=\Lin{(x,x_k^*)x_k}_{k \in N_1}$, since $(x,x_k^*)\neq 0$ for $k \in N_1$. 
\end{proof}

\begin{Example} \label{e_1+e_k}
 Let  $\{e_k\}_{k=1}^\infty$ be an orthonormal basis. Let $x_k=e_1+e_k$ for $k\geq2$. The system $\{x_k\}_{k=2}^\infty$ is complete and minimal but not hereditarily complete.

 Indeed, consider an element $x$ that is orthogonal to all $x_k$. Let us represent it as  $x=\sum\limits_{k=1}^\infty c_ke_k$. From orthogonality, we obtain the system of equations $c_1+c_m=0$, $m \geq 2$, which implies $-c_1=c_2=\ldots=c_m=\ldots$ The coefficients of the Fourier series are square-summable, which means they tend to zero, therefore, $c_1=0$, and so $c_k=0 $ for all $k\in \mathbb{N}$. The latter proves the required completeness.  

 Obviously, for $\{x_k\}_{k=2}^\infty$, the biorthogonal sequence will be $\{e_k\}_{k\geq 2}$. The system $\{e_k\}_{k\geq 2}$ is incomplete, and then the system $\{x_k\}_{k=2}^\infty$ is not hereditarily complete.
\end{Example}
For the sake of completeness and to illustrate the methods used below, we present a proof of Young's result \cite{RY}.  
\begin{Lemma} \label{example_inf}
Let $H_0$ be an infinite-dimensional closed subspace in $H$. Then there exists a complete minimal system $\{x_k\}_{k=1}^\infty$ for which $\overline{Lin} \{x_k^* \}_{k=1}^\infty =H_0$.
\end{Lemma}
\begin{proof}
Let us use the representation $H=H_0 \oplus H_0^\bot$. Let us choose some orthonormal bases $\{e_k\}_{k=1}^\infty$ in $H_0$ and $\{f_m\}_{m=1}^\infty$ in $H_0^{\bot}$ (if $H_0^{\bot}$ is finite-dimensional, then first we take vectors from the basis as $f_m$, and fill in the remaining places with zeros). Now we introduce a vector system, assuming $x_k=2^k\sum \limits _{j=1}^k\frac{1}{k^{j-1}}f_j + e_k$, where $k\in \mathbb{N}$. We shall prove that this system is complete. Any element $x$ from $H$ can be represented as $x=\sum \limits_{k} c_k e_k + \sum \limits_{m} d_m f_m $. Here we assume that $d_m=0$ if $f_m=0$. Let the vector $x$ be orthogonal to all vectors $x_k$. Then 
$$
2^k\sum \limits_{j=1}^k\frac{d_j}{k^{j-1}}+c_k=0 \Longleftrightarrow d_1+\frac{1}{k}\sum \limits_{j=2}^k\frac{d_j}{k^{j-2}}+\frac{c_k}{2^k}=0.
$$
 Moreover, $\{ c_k\}_{k=1}^\infty \in \ell^2$, so $c_k\xrightarrow[k\rightarrow \infty]{}0$; $\{ d_m\}_{m=1}^\infty \in \ell^2$, from which we conclude that the inequality $| d_m| \leq M$ is true for all $m$ with some constant $M$. Hence
$$
\left| \frac{1}{k}\sum \limits_{j=2}^k\frac{d_j}{k^{j-2}}\right| \leq \frac{M}{k}\sum \limits_{j=2}^k\frac{1}{k^{j-2}} \leq \frac{M}{k} \frac{1}{1-\frac{1}{k}}=\frac{M}{k-1} \xrightarrow[k \rightarrow \infty]{} 0.
$$
Therefore, $d_1=0$. From the equality $2d_1 + c_1=0$, we get that $c_1=0$. Continuing these arguments by induction, we see that all coefficients in the expansion of $x$ are zero. 

We will prove minimality by presenting a biorthogonal system: as can be easily verified, the system $\{e_k\}_{k=1}^\infty$ fits. By definition, $\overline{Lin}\{e_k\}_{k=1}^\infty = H_0$.
\end{proof}

\section{Vector replacement property and criterion of hereditary completeness}

Let us prove a useful lemma that will allow us to find a criterion for hereditary completeness based only on the vector system, but not on its biorthogonal one. 
\begin{Lemma} \label{swap}
Let $\{x_k\}_{k=1}^\infty$ be a system of vectors having biorthogonal $\{x_k^*\}_{k=1}^\infty$, $N_1 $ be a subset of $\mathbb{N}$, $N_2=N_1^c$, and $k_0 \in N_2$. Let us also assume that the defect of the system $\{x_k\}_{k=1}^\infty$ with respect to $N_1$ is equal to $m \in \Z_+ \cup \{\infty\}$. Then if we replace $N_1$ with $N_1\cup \{k_0\}$, and $N_2$ with $N_2 \setminus \{k_0\}$, the system defect will not change. In particular, if the system $\{x_k\}_{k \in N_1}\cup \{x_k^*\}_{k \in N_2}$ was complete, it will remain so after replacement.
\end{Lemma}
\begin{proof}
The Hilbert space $H$ can be represented as
\begin{equation*}
    H=\overline{Lin} \{x_k\}_{k \in N_1}\oplus \overline{Lin} \{x_k^*\}_{k \in N_2}\oplus W,
\end{equation*}
where $W$ is a subspace of $H$ having dimension $m$. Let us write another representation:
\begin{equation*}
    H=\overline{Lin} \{x_k\}_{k \in N_1}\oplus \overline{Lin} \{x_k^*\}_{k \in N_2 \setminus \{k_0\}}\oplus W'.
\end{equation*}
We shall show that $\dim W'=m+1$ (we assume that $\infty +1=\infty$). Indeed, the vector $x_{k_0}^*$ decomposes as $x_{k_0}^*=u+v$, where $u \in \overline{Lin} \{x_k^*\}_{k \in N_2 \setminus \{k_0\}}$, and $v $ is orthogonal to $ \overline{Lin} \{x_k^*\}_{k \in N_2 \setminus \{k_0\}}$. The vector $v$ is not equal to $0$, otherwise $(x_{k_0},x_{k_0}^*)=(x_{k_0},u)=0 \neq 1$. In addition, if the vector $h$ lies in the subspace $\overline{Lin} \{x_k^*\}_{k \in N_2}$ and is orthogonal to the subspace $ \overline{Lin} \{x_k^*\}_{k \in N_2 \setminus\{k_0\}}$, as well as the vector $v$, then $h$ is orthogonal to $ x_{k_0}^*$, from which we get $h=0$. Therefore, we have a representation
\begin{equation*} 
 \overline{Lin} \{x_k^*\}_{k \in N_2}=\overline{Lin} \{x_k^*\}_{k \in N_2 \setminus\{k_0\}} \oplus Lin\{v\}.   
\end{equation*}
Hence $W'=Lin\{v\}\oplus W$, which gives the required identity for the dimension of $W'$. Now we write the identity $x_{k_0}=a+b+c+d$ corresponding to the decomposition
\begin{equation*} 
   H=\overline{Lin} \{x_k\}_{k \in N_1}\oplus \overline{Lin} \{x_k^*\}_{k \in N_2 \setminus \{k_0\}}\oplus Lin\{v\} \oplus W.              
\end{equation*}
The vector $b$ is equal to $0$ because $x_{k_0}$ is orthogonal to $ \overline{Lin} \{x_k^*\}_{k \in N_2 \setminus \{k_0\}}$, and $c\neq 0$, since
\begin{equation*}
  1=(x_{k_0},x_{k_0}^*)=(a+c+d,u+v)=(a+c+d,v)=(c,v).  
\end{equation*}
Let us write one more decomposition of $H$ into a direct sum:
\begin{equation*}
 H=\overline{Lin} \{x_k\}_{k \in N_1\cup \{k_0\}}\oplus \overline{Lin} \{x_k^*\}_{k \in N_2 \setminus \{k_0\}}\oplus W''.   
\end{equation*}
It is easy to see that $W''=\{h\in W' \ |\  h \bot x_{k_0} \}$. But any $h$ from $ W'$ is orthogonal to $a$, therefore $W''=W' \ominus Lin\{c+d\}$, moreover, $c+d\neq 0$. The last equality gives $\dim W''=\dim W' -1=m$, which is what was required.
\end{proof}

\begin{Remark}
   In Lemma \ref{swap} the roles of the elements of the original system and its biorthogonal can be reversed, since one of the biorthogonal to $x_k^*$ is the system $x_k$.
\end{Remark}
\begin{Corollary} \label{swap_fin}
If we change the partition $ \mathbb{N}=N_1 \sqcup N_2$ by a finite number of elements (that is, exclude some finite $N \subset N_1$ from $N_1$ and add it to $N_2$ or vice versa), then the completeness or incompleteness of the system $\{x_k\}_{k \in N_1} \cup \{x_k^*\}_{k \in N_2} $ will not change, moreover, the dimension of the orthogonal complement to the closed linear span of the system will not change.
\end{Corollary}
We will need the following standard facts about the convergence of monotone sequences of projectors.
\begin{Lemma}
\label{mon_conv}
\begin{enumerate}
    \item  \label{mon_conv_1}  Let $\{ x_k\}_{k \in I_j}$, $j=1, 2 \ldots$, be sets of vectors from $H$, with $I_j \subset I_{j+1}$, 
$$
V_j=\overline{Lin} \{ x_k\}_{k \in I_j}, \ V=\overline{Lin} \{ x_k\}_{k \in \bigcup \limits_{j=1}^\infty I_j}.
$$
Then $P_{V_j} \xrightarrow[m \rightarrow \infty]{} P_V$ pointwise.  
\item \label{mon_conv_2} Let $\{V_j\}_{j=1}^\infty$ be closed subspaces, $V_{j+1} \subset V_j$, and
$\bigcap\limits_{j=1}^\infty V_j=V $. Then $ P_{V_j} \xrightarrow[j \rightarrow \infty]{}P_V$ pointwise. 
\end{enumerate}

\end{Lemma}

We shall introduce some notation that will be used in what follows. For $\sigma\subset \mathbb{N}$, we put
\begin{equation} \label{sigma_m}
    \sigma_m:=\sigma \cup [m+1:+\infty)
\end{equation}
 and, consequently,
\begin{equation} \label{sigma_m^c}
    (\sigma_m)^c=\sigma^c \cap [1:m]. 
\end{equation}
Here $[a:b]$ denotes all integers $n$ satisfying the inequalities $a\leq n\leq b$, and $[a:+\infty)$ --- all integers $n$ with the condition $n\geq a$.

For a minimal system $\{x_k\}_{k=1}^\infty$ and $\sigma \subset \mathbb{N}$ we put $H_\sigma=\Lin{x_k}_{k \in \sigma}$ and $H'_\sigma=\Lin{x_k^*}_{k \in \sigma^c}$. In addition, we define $P_\sigma=P_{H_\sigma}$ and $P'_\sigma=P_{H'_\sigma}$. By definition of biorthogonal system, subspaces $H_\sigma$ and $H'_\sigma$ are orthogonal, which means $P_\sigma + P'_\sigma$ is the projector on $H_\sigma\oplus H'_\sigma$. 
\begin{Theorem} \label{mixed_sys}
Let $\{x_k\}_{k=1}^\infty$ be a complete minimal system, and $\sigma$ be a subset of the natural numbers. Then $\{x_k\}_{k \in \sigma } \cup \{ x_k^*\}_{k \in \sigma^c}$ is complete if and only if 
$$
\bigcap\limits_{m=1}^\infty H_{\sigma_m}=H_{\sigma}.
$$
\end{Theorem}
\begin{proof}
 Let us prove that the completeness of the mixed system $\{x_k\}_{k \in \sigma } \cup \{ x_k^*\}_{k \in \sigma^c}$ is equivalent to pointwise convergence of projectors $P_{\sigma_m} \xrightarrow[m \rightarrow \infty]{}P_\sigma$.

First, let us assume that the convergence of projectors takes place. The subset $\sigma_m$  differs from $\mathbb{N}$ by a finite number of elements, so by Corollary \ref{swap_fin} mixed system $\{x_k\}_{k\in \sigma_m} \cup \{x_k^*\}_{k\in (\sigma_m)^c}$ is complete. The latter is equivalent to the equality $P_{\sigma_m}+P_{\sigma_m}'=I$. The inclusions $(\sigma_{m+1})^c \supset (\sigma_{m})^c$ and the equality $\bigcup\limits_{m=1}^\infty (\sigma_m)^c=\sigma^c$ hold. Therefore from the statement \ref{mon_conv_1} of Lemma \ref{mon_conv} we get $P_{\sigma_m}' \xrightarrow[m \rightarrow \infty]{} P'_\sigma$. From here we conclude, passing in the equality above to the limit, that $P_\sigma + P'_\sigma=I$, and hence $H_\sigma \oplus H'_\sigma=H$, that is, the mixed system is complete.

Let us assume that the completeness of the system  $\{x_k\}_{k \in \sigma } \cup \{ x_k^*\}_{k \in \sigma^c}$ is known. As above, $P_{\sigma_m}+P'_{\sigma_m}=I$ and $P'_{\sigma_m} \xrightarrow[m \rightarrow \infty]{} P'_\sigma$. By the previous reasoning $P_{\sigma_m}=I -P'_{\sigma_m} \xrightarrow[m \rightarrow \infty]{}I - P'_\sigma$. Due to the completeness of the mixed system, $I-P'_\sigma=P_\sigma$. From here we get the necessary convergence. 

Now the result follows easily from the statement \ref{mon_conv_2} of Lemma \ref{mon_conv}, because $\sigma_{m+1} \subset \sigma_m$ and $P_{\sigma_m} \xrightarrow[m \rightarrow \infty]{} P$, where $P$ is the projector on $\bigcap\limits_{m=1}^\infty H_{\sigma_m}$.
\end{proof}

\begin{Corollary} \label{hered_comp_proj}
The system $\{x_k\}_{k=1}^\infty$ is hereditarily complete if and only if the equalities
$$
\bigcap\limits_{m=1}^\infty H_{\sigma_m}=H_{\sigma}
$$ 
hold for any subsets $\sigma$ of natural numbers, which, in turn, is equivalent to pointwise convergence
$$
P_{\sigma_m} \xrightarrow[m \rightarrow \infty]{} P_\sigma.
$$
\end{Corollary}

Another, topological, characterization of hereditarily complete systems is possible. Let us define a topology on all possible subsets $\sigma \subset \mathbb{N}$. To each subset $\sigma$ the sequence $\{a_k^{\sigma}\}_{k=1}^\infty$ of zeros and ones can be assigned according to the following rule: $a_k^{\sigma}:=1$ if $k \in \sigma$ and $a_k^{\sigma}:=0$ otherwise. It is clear that the mapping $\sigma \xrightarrow[]{} a^{\sigma}$ defines a one-to-one correspondence between the set $\mathfrak{P}(\N)$ of all subsets of $\mathbb{N}$ and $\{0,1\}^{\mathbb{N}}$. Now we define a metric on the set of sequences of zeros and ones. 
$$
\rho(a,b):=\sum \limits_{k=1}^\infty \frac{|a_k - b_k|}{2^k}.
$$
It is known that this metric determines on a set of sequences (and hence on a set of subsets of natural numbers) the structure of a compact topological space with the topology of the product of a countable number of two-point discrete spaces. It is also known that the convergence $a^j \xrightarrow[j \rightarrow \infty]{}$ $a$ is equivalent to the following property: $\forall m \in \mathbb{N} \ \  \exists N=N(m): \forall j\geq N \ \  a_i^{j}=a_i$ for $1\leq i \leq m$ (in terms of subsets, this means that starting from a certain number they match on the first $m$ natural numbers). It is easy to show that the operations of intersection, union, and complementation are continuous on this space as functions of two variables.

Let $\{x_k\}_{k=1}^\infty$ be a complete minimal system.  Let us also assume that the norms of all $x_k$ are equal to 1. Normalization will not change the linear spans of the vectors $x_k$ or $x_k^*$ in any way, but it will be convenient later. Let us map each subset $\sigma \subset \mathbb{N}$  to the projector $P_\sigma$ on $\overline{Lin}\{x_k\}_{k \in \sigma}$. It is clear that the correspondence $\sigma \mapsto P_\sigma$ is injective. Indeed, let $P_{\sigma^1}=P_{\sigma^2}$ and, for instance, $k\in \sigma^1$, but $k \notin \sigma^2$. In this case, $x_k \in \overline{Lin}\{x_k\}_{k \in \sigma^2}$, which contradicts the minimality of the system. Let us define a topology induced by a strong operator topology on the set of all orthogonal projectors onto closed subspaces of $H$. It follows from the Banach--Steinhaus theorem that this topology is also metrizable and can be given by a metric
$$
d_s(P_1,P_2)=\sum \limits_{k=1}^\infty \frac{\|P_1x_k - P_2 x_k \|}{2^k}.
$$
In addition, we will consider the set of projectors endowed with a weak operator topology. It is clear that this topology is generated by the metric
$$
d_w(P_1,P_2)=\sum \limits_{k,j=1}^\infty \frac{|(P_1 x_k - P_2 x_k ,x_j)|}{2^{k+j}}.
$$
\begin{Lemma}
\label{hered_comp_top}
 The set of all projectors $P_\sigma$ corresponding to a complete minimal system $\{x_k\}_{k=1}^\infty$ will be written as $\Omega(x_k)$. Let us denote the previously defined mapping $\sigma \mapsto P_{\sigma}$ from $\mathfrak{P}(\mathbb{N})$ to $\Omega (x_k)$ as $\phi$.   
\begin{enumerate}
    \item  \label{hered_comp_top_1} The mapping $\phi^{-1}: \Omega(x_k) \rightarrow  \mathfrak{P}(\mathbb{N})$ is uniformly continuous on the set $\Omega(x_k)$ equipped with the metric corresponding to the strong operator topology or the metric corresponding to the weak operator topology. In addition, if $\sigma$ is such that the system $\{x_k\}_{k \in \sigma}\cup \{x_k^*\}_{k \in \sigma^c}$ is complete in $H$, then $\phi$ is continuous at the point $\sigma$ (projectors can be equipped with either of the metrics).
    \item \label{hered_comp_top_2} The sets $\sigma$ corresponding to complete mixed systems are exactly the points of continuity of the mapping $\phi$. In particular, the mapping $\phi$ is a homeomorphism if and only if the system $\{x_k\}_{k=1}^\infty$ is hereditarily complete. 
    \item The system $\{x_k\}_{k=1}^\infty$ is hereditarily complete if and only if the projector space $\Omega(x_k)$ is compact.
\end{enumerate} 

\end{Lemma}
\begin{Remark}
    We put an emphasis on the notation. Here and below, a notation of the form $\sigma^m$ with an index at the top denotes some sequence of sets (chosen in accordance with a certain property), and the notation $\sigma_m$ denotes a specific sequence of sets constructed from one set $\sigma$, as was done in the equality \eqref{sigma_m}.  
\end{Remark}

\begin{proof}
    (1) Note that 
\begin{multline*}
  d_w(P_1,P_2)=\sum \limits_{k,j=1}^\infty \frac{|(P_1 x_k - P_2 x_k ,x_j)|}{2^{k+j}}\leq \sum \limits_{k,j=1}^\infty \frac{\|P_1 x_k - P_2 x_k\| \cdot \|x_j\|}{2^{k+j}}= \\ =\sum \limits_{k=1}^\infty \frac{\|P_1 x_k - P_2 x_k\| }{2^{k}}  \cdot \left(\sum \limits_{j=1}^\infty \frac{1}{2^j}\right) =d_s(P_1,P_2).  
\end{multline*}
Therefore, it is sufficient to prove the uniform continuity of $\phi^{-1}$ only with respect to the metric $d_w$. Let us assume the opposite. Let there be a sequence of pairs of projectors $P_{\sigma^m}$, $P_{\tau^m}$ and such $\varepsilon >0$ that $d_w(P_{\sigma^m}, P_{\tau^m}) \xrightarrow[m \rightarrow \infty]{} 0$, but $\rho(\sigma^m,\tau^m)\geq \varepsilon $. In this case, there exists a number $m_0$ such that for any $m \in \mathbb{N}$ the set $\sigma^{m} $ does not coincide with $\tau^{m}$ on the segment of natural numbers from 1 to $m_0$. Without loss of generality, we can assume that there is $p \in \sigma^{m} $ but $p \notin \tau^m $ for all $m$. Indeed, since there are infinitely many different $m$, and finitely many numbers from 1 to $m_0$, there is a subsequence in $\mathbb{N}$ such that $\sigma^{m}$ and $\tau^m$ differ on this subsequence by a number with index $1\leq p \leq m_0$. We consider only the case $p \in \sigma^{m} $, $p \notin \tau^m$, the second case is considered similarly. We have 
\begin{align*}
  d_{w}(P_{\sigma^m}, P_{\tau^m}) & \geq \frac{|(P_{\sigma^{m}}x_p - P_{\tau^m} x_p,x_p)|}{2^{2p}}=\frac{|(x_p - P_{\tau^m} x_p,x_p)|}{2^{2p}}= \\ &=\frac{\|x_p - P_{\tau^m} x_p\|^2}{2^{2p}}\geq  \frac{\|x_p - P_{\mathbb{N} \setminus \{p\} }  x_p\|^2}{2^{2p}}>0,    
\end{align*}
  which gives a contradiction. 

  Let us prove the second part of (1). We check that the projectors $P_{\sigma^m}$ converge pointwise to $P_\sigma$ if $\sigma^m \rightarrow \sigma$ (and then they converge in the weak topology as well). It suffices to check this only for vectors from the complete system, then the convergence will also hold on their linear span, and then for any vector from $H$ by the Banach--Steinhaus theorem. As a complete system we take $\{x_k\}_{k \in \sigma}\cup \{x_k^*\}_{k \in \sigma^c}$. Let $n \in \sigma$, then for $m$ sufficiently large the sets  $\sigma^m$ will contain $n$, but then $P_{\sigma^m}x_n=x_n=P_{\sigma}x_n$. Similarly, for any number $n \in  \sigma^c$, for sufficiently large $m$ the sets $\sigma^m$ do not contain $n$, which implies $P_{\sigma^m}x_n^*=0=P_{\sigma}x_n^*$. This completes the proof.
  
  (2) Let $\sigma$ be such that the system $\{x_k\}_{k \in \sigma}\cup \{x_k^*\}_{k \in \sigma^c}$ is not complete. Consider the family of sets $\sigma_m$. They coincide with $\sigma$ at least on the first $m$ natural numbers, and hence tend to $\sigma$ in the topology under consideration. On the other hand, we know from the statement \ref{mon_conv_2} of the Lemma \ref{mon_conv} on the limit of decreasing projections that $P_{\sigma_m} \xrightarrow[m \rightarrow \infty]{} P$, where $P$ is the projection onto the intersection of the spaces $\overline{Lin}\{x_k\}_{k \in \sigma_m}$, and this intersection is not equal to $\overline{Lin}\{x_k\}_{k \in \sigma}$ in the case when the mixed system corresponding to $\sigma$ is not complete by Theorem \ref{mixed_sys}. 

  (3) In one direction the statement has already been proven, in the other direction it follows immediately from the fact that $\phi^{-1}$ is a continuous bijection from a compact to a Hausdorff space. 
\end{proof}

In the general case, the space of projections $P_\sigma$ of a complete minimal system is not homeomorphic to the space of subsets of natural numbers, but the topology of the space of projections $P_\sigma$ can be represented using the graph of a bounded function on a countable product of the two-point spaces.  
\begin{Lemma} \label{graph_proj}
The following statements are equivalent:
\begin{enumerate}
    \item Convergence $P_{\sigma^m} \xrightarrow[m \rightarrow \infty]{} P_{\sigma}$ takes place. 
    \item Two limit relations are true: $\sigma^{m} \xrightarrow[m \rightarrow \infty]{} \sigma$ and $d_s(P_{\sigma^m},0) \xrightarrow[m \rightarrow \infty]{} d_s(P_{\sigma},0)$.
\end{enumerate}
  The convergence of projectors here can be understood both in the sense of weak and in the sense of strong topology. 
\end{Lemma}
\begin{proof}
    (1) Note, first, that it is sufficient to prove the criterion only for pointwise convergence. If the projectors converge pointwise, then they also converge weakly. On the other hand, if the weak limit of the projectors is a projector, then pointwise convergence also holds (see, for example, \cite{BS}, p. 54, Theorem 8). If $P_{\sigma^m} \xrightarrow[m \rightarrow \infty]{} P_{\sigma}$, then from Lemma \ref{hered_comp_top} we know that $\sigma^{m} \xrightarrow[m \rightarrow \infty]{} \sigma$, and $d_s(P_{\sigma^m},0) \xrightarrow[m \rightarrow \infty]{} d_s(P_{\sigma},0)$ from the general properties of convergence in metric spaces.
    
    (2) We will prove the converse by contradiction. Let us assume that $P_{\sigma^m}$ do not converge to $P_{\sigma}$, then there exist $\varepsilon >0$ and $j_0 \in \mathbb{N}$ such that $\|P_{\sigma^m}x_{j_0} - P_{\sigma}x_{j_0} \|\geq \varepsilon $ for an infinite number of indices $m$. Passing to a subsequence, we will assume that the inequality is satisfied for any natural $m$. Let $N(m)$ be the largest number such that $\sigma \cap [1:N(m)]=\sigma^m \cap [1:N(m)]$ (if the equality is true for all $N$, then we set $N(m)=m$). From the convergence of $\sigma^m$ to $\sigma$ it follows that $N(m)  \xrightarrow[m \rightarrow \infty]{} \infty$. For any $j$
    $$
    \|P_{\sigma^m} x_j\| \geq \|P_{\sigma^m \cap [1:N(m)]} x_j\|=\|P_{\sigma \cap [1:N(m)]}x_j\|\xrightarrow[m \rightarrow \infty]{} \|P_{\sigma}x_j\|
    $$
    (see assertion \ref{mon_conv_1} of Lemma \ref{mon_conv}). From this we obtain the estimate 
    \begin{equation}
         d_s(P_{\sigma^m},0)= \sum \limits_{j=1}^\infty \frac{\|P_{\sigma^m}x_j\|}{2^j}\geq \frac{\|P_{\sigma^m} x_{j_0}\|}{2^{j_0}} + \sum \limits_{j=1, j\neq j_0}^{\infty}\frac{\|P_{\sigma \cap [1:N(m)]}x_j\|}{2^j},  \label{d_s(P_sigma^m,0)}
    \end{equation}
    moreover, the second term on the right-hand side tends to $\sum \limits_{j=1, j\neq j_0}^{\infty}\frac{\|P_{\sigma}x_j\|}{2^j}$ as $m$ tends to infinity, since the sum is majorized by the series $\frac{1}{2^j}$ regardless of $m$. Let us now estimate the first term on the right-hand side. We write the identities  
    $$
    P_{\sigma}x_{j_0}=P_{\sigma \cap [1:N(m)]}x_{j_0} + (P_{\sigma} - P_{\sigma \cap [1:N(m)]})x_{j_0}
    $$
    and
    $$
    P_{\sigma^m}x_{j_0}=P_{\sigma^m \cap [1:N(m)]}x_{j_0} + (P_{\sigma^m} - P_{\sigma^m \cap [1:N(m)]})x_{j_0}=P_{\sigma \cap [1:N(m)]}x_{j_0} + (P_{\sigma^m} - P_{\sigma^m \cap [1:N(m)]})x_{j_0}.
    $$
    From these equalities we obtain 
    \begin{multline*}
       \varepsilon  \leq\|P_{\sigma^m}x_{j_0} - P_{\sigma}x_{j_0}\|=\|(P_{\sigma^m} - P_{\sigma^m \cap [1:N(m)]})x_{j_0} - (P_{\sigma} - P_{\sigma \cap [1:N(m)]})x_{j_0}\|\leq \\ \leq  \|(P_{\sigma^m} - P_{\sigma^m \cap [1:N(m)]})x_{j_0}\|+\|(P_{\sigma} - P_{\sigma \cap [1:N(m)]})x_{j_0}\|.
    \end{multline*}
 But according to the statement \ref{mon_conv_1} of the Lemma \ref{mon_conv} we get $\|(P_{\sigma} - P_{\sigma \cap [1:N(m)]})x_{j_0}\| \xrightarrow[m \rightarrow \infty]{} 0 $, whence we find $\lowlim \limits_{m \to \infty} \|(P_{\sigma^m} - P_{\sigma^m \cap [1:N(m)]})x_{j_0}\| \geq \varepsilon $. Considering now that
 $$
 \|P_{\sigma^m}x_{j_0}\|^2=\|P_{\sigma \cap [1:N(m)]}x_{j_0}\|^2+\|(P_{\sigma^m} - P_{\sigma^m \cap [1:N(m)]})x_{j_0}\|^2,
 $$
 we conclude 
 $$
 \lowlim \limits_{m \to \infty} \|P_{\sigma^m}x_{j_0}\|\geq \sqrt{\|P_{\sigma}x_{j_0}\|^2+\varepsilon ^2}>\|P_{\sigma}x_{j_0}\|.
 $$
 From the last inequality and from the estimate \eqref{d_s(P_sigma^m,0)} it follows that $\lowlim \limits_{m \to \infty}  d_s(P_{\sigma^m},0)>d_s(P_{\sigma},0)$ --- a contradiction. 
 \end{proof}
 \begin{Remark}
      From the proof of the Lemma \ref{graph_proj} it is clear that if the convergence $P_{\sigma^m} \xrightarrow[m \rightarrow \infty]{} P_{\sigma}$ does not hold, then $\lowlim \limits_{m \to \infty}  d_s(P_{\sigma^m},0)>d_s(P_{\sigma},0)$, that is, the function $\sigma \mapsto d_s(P_\sigma,0)$ is lower semicontinuous on the space of subsets $\mathfrak{P}(\mathbb{N})$. 
 \end{Remark}
 
As an example of the application of the developed technique, we will prove another theorem on hereditary completeness.
\begin{Theorem} \label{mix_hered_comp}
Let $\{x_k\}_{k=1}^\infty$ be a hereditarily complete system. Then for any $ \sigma \subset \mathbb{N}$ the system $\{x_k\}_{k \in \sigma } \cup \{ x_k^*\}_{k \in \sigma^c}$ is also hereditarily complete.
\end{Theorem}
\begin{proof}
From the hereditary completeness of the system follows the hereditary completeness of its biorthogonal (the criterion of hereditary completeness from the Lemma \ref{hered_comp_crit} is formulated symmetrically). Let us denote
\begin{equation*}
w_k= 
 \begin{cases}
   x_k, &\text{$k\in \sigma$},\\
      x_k^*, &\text{ $k\in \sigma^c $}.
 \end{cases}
\end{equation*}
Let $\tau,\tau^m \subset \mathbb{N}$ and $\tau^m \xrightarrow[m \rightarrow \infty]{} \tau$. We introduce the notations $\tau^1:=\tau \cap \sigma$, $\tau^2:=\tau \cap \sigma^c$, and also $\tau^{1,m}:=\tau^m \cap \sigma $, $\tau^{2,m}:=\tau^m \cap \sigma^c$. The system $\{w_k\}_{k=1}^\infty$ is complete and minimal. Completeness follows from the hereditary completeness of $\{x_k\}_{k=1}^\infty$. Minimality is obtained if we note that the minimality of the systems $\{x_k\}$ and $\{x_j'\}$, as well as the orthogonality of all $x_k$ to all $x_j'$, imply the minimality of the system $\{x_k\}\cup \{x_j'\}$. Indeed, the vectors biorthogonal to the system can be taken from the closed linear span of the system itself. For example, we can consider the normalized vectors $P_{\Lin{x_l}_{l \neq k}} x_k$. 

Let us introduce one new notation. If it is not clear what minimal system $\{y_k\}$ we are talking about, we will write $H_{\sigma, y_k}=\Lin{y_k}_{k \in \sigma}$ and $P_{\sigma,y_k}=P_{H_{\sigma, y_k}}$. From the orthogonality of $x_k$ and $x_j^*$ for $k \neq j$ we get the representations 
$$
H_{\tau, w_k}=H_{\tau^1, w_k} \oplus H_{\tau^2,w_k}=H_{\tau^1, x_k}\oplus H_{\tau^2,x_k^*}
$$
 and
 $$
 H_{\tau^m, w_k}=H_{\tau^{1,m}, x_k}\oplus H_{\tau^{2,m},x_k^*}.
 $$
 Moreover, $P_{\tau^{1,m},x_k} \xrightarrow[m \rightarrow \infty]{} P_{\tau^1,x_k}$ and  $P_{\tau^{2,m},x_k} \xrightarrow[m \rightarrow \infty]{} P_{\tau^2,x_k^*}$ pointwise due to hereditary completeness of the systems $\{x_k\}_{k=1}^\infty$ and $\{x_k^*\}_{k=1}^\infty$ by the Lemma \ref{hered_comp_top}, as well as by the continuity of the intersection operation in the topology on subsets of natural numbers. Obviously, if $Z_m= X_m \oplus Y_m$ is a sequence of subspaces, and there is pointwise convergence $P_{X_m} \xrightarrow[m \rightarrow \infty]{} P_X$, $P_{Y_m} \xrightarrow[m \rightarrow \infty]{} P_Y$, then $P_{Z_m} \xrightarrow[m \rightarrow \infty]{} P_{X \oplus Y}$. But then $P_{\tau^m,w_k} \xrightarrow[m \rightarrow \infty]{} P_{\tau, w_k}=P_{\tau^1, x_k}+P_{\tau^2,x_k^*}$. Since the sequence $\tau^m \xrightarrow[m \rightarrow \infty]{} \tau$ was arbitrary, we complete the proof by applying Lemma \ref{hered_comp_top}, statement 2. 
\end{proof}

\section{Systems with prescribed sets of defects}
In view of the Lemma \ref{hered_comp_crit}, the question arises about what defects a complete minimal system can have (see Definition \ref{def}). It turns out that its set of defects relative to all possible subsets $\sigma$ can be almost arbitrary.

\begin{Theorem} \label{def_custom}
    Let $S=\{0=k_0, k_1, \ldots, k_s, \ldots,\}$ be either a finite set of non-negative integers or a set containing infinity (either finite or infinite). Then there exists a complete minimal system whose defects with respect to all possible subsets of the natural numbers form exactly the set $S$.
\end{Theorem}
\begin{proof}
We shall divide the proof into several steps. \medskip

\noindent  \textbf{Step 1: case of defects 0, m} \medskip

First, we show that there exist complete minimal systems of vectors such that the set of their defects consists of only two numbers 0 and $m$. Consider the system of vectors
$$
x_k=e_1 + ke_2 +\ldots +k^{m-1}e_m + e_{m+k}.
$$
Systems of this type generalize the Example \ref{e_1+e_k}. It is easy to see that its biorthogonal is the system $\{e_{m+k}\}_{k=1}^\infty$. Let us prove the completeness of the system $\{x_k\}_{k=1}^\infty$. Let $x=\sum \limits_{j=1}^\infty c_j e_j$ be orthogonal to all $x_k$. Then $\{\sum \limits_{j=1}^m c_j k^{j-1}\}_{k=1}^\infty \in \ell^2$ must hold (since this is the expression for the coefficients $-c_{m+k}$). In particular, $\sum \limits_{j=1}^m c_j k^{j-1} \xrightarrow[k \rightarrow \infty]{} 0$. Dividing each of the sums by $k^{m-1}$, we obtain $\sum \limits_{j=1}^{m-1} c_j k^{j-m}+c_m \xrightarrow[k \rightarrow \infty]{} 0$, whence $c_m=0$. Knowing that $c_m=0$, we divide the sum $\sum \limits_{j=1}^{m-1} c_j k^{j-1}$ by $k^{m-2}$ and repeat the reasoning. Continuing the process, we obtain that $c_j$ equals 0 for $j$ from 1 to $m$, and then $x$ equals 0, because all its Fourier coefficients are expressed through $c_1, \ldots, c_m$.

Consider the system $\{x_k\}_{k \in \sigma } \cup \{x_k^*\}_{k \in \sigma^c}$. If $\sigma$ is finite, then the corresponding defect is equal to $m$ by Corollary \ref{swap_fin}. If $\sigma$ is infinite, then for the vector $x$ orthogonal to the system the coefficient $c_{m+k}$ equals 0 for $k \in \sigma^c$, since $x_k^*=e_{m+k}$, and, having carried out the same reasoning for $\{x_k\}_{k \in \sigma}$ as above, we see that all other coefficients in the expansion of $x$ are equal to zero. Thus, the defect in this case is equal to zero. \medskip

\noindent  \textbf{Step 2: finite set of integer defects} \medskip

Now let us construct a complete minimal system with an arbitrary finite set of defects consisting of non-negative integers $\{0=k_0, k_1, k_2, \ldots, k_s\}$ (the numbers in the set are arranged in ascending order). We set
$$
x_k^j=k^{k_j}e_{k_j +1}+k^{k_j+1}e_{k_j +2}+\ldots+k^{k_s -1}e_{k_s}+e_{k+k_s}, \quad j=0,\ldots, s-1,
$$
$$
x_k^s=e_{k+k_s},
$$
$k=1, 2 \ldots$ Now we denote $x_k=x_k^j$ for $k-1 \equiv j $ (mod $s+1$), that is, the sequence $x_k$ includes the following vectors $x_1^0, x_2^1, \ldots, x_{s+1}^s, x_{s+2}^0, \ldots$ The sequence $x_k$ has biorthogonal $\{e_{k+k_s}\}_{k=1}^\infty$ and is complete. To establish its completeness, it is sufficient to consider $x$ orthogonal to all vectors $\{x_{(s+1)n +1 }\}_{n=0}^\infty$, deduce from here by the same method as before that all Fourier coefficients for $e_1, \ldots, e_{k_s}$ are equal to zero, and from this obtain that any $x$ orthogonal to all $x_k$ is equal to zero. 

Let us prove that the obtained sequence has the required set of defects. Let $\sigma \subset \mathbb{N}$, $\sigma=\sigma_{j_1}\sqcup \ldots \sqcup \sigma_{j_r} $, where $\sigma_{j_p}=\sigma \cap \{(s+1)n+j_p+1\}_{n=1}^\infty$ are the sets on which $x_k$ equal $x_k^j$ with fixed $j$, $j_1<\ldots<j_r$. Without loss of generality, we assume that all $\sigma_{j_p}$ are infinite, otherwise we exclude $\sigma_{j_p}$ from $\sigma$ and add to $\sigma^c$. The defect will not change from this by the Corollary \ref{swap_fin}.  Consider the vector $x=\sum \limits_{j=1}^\infty c_j e_j$, orthogonal to the system $\{x_k\}_{k \in \sigma} \cup \{x_k^*\}_{k \in \sigma^c}$. This vector $x$ is orthogonal to an infinite number of vectors $x_k^{j_1}$ (namely, $x_k$, $k \in \sigma_{j_1}$), and therefore there are an infinite number of equations $\sum \limits_{l=k_{j_1+1}}^{k_s} c_l k^{l-1} + c_{k+k_s}=0$. From them we obtain $c_{k_{j_1}+1}=c_{k_{j_1}+2}=\ldots=c_{k_s}=0$ and $c_{k+k_s}=0$ for $k \in \sigma_{j_1}$. The orthogonality condition of $x$ and $x_k$, $k \in \sigma_{j_p}$, where $p>1$, yields the equations $\sum \limits_{l=k_{j_p+1}}^{k_s} c_l k^{l-1} + c_{k+k_s}=0$. Moreover, $k_{j_p}>k_{j_1}$, which means $c_{k+k_s}=0$ also for $k \in \sigma_{j_p}$. Finally, $c_{k+k_s}=0$ for $k \in \sigma^c$, since $x_k^*=e_{k+k_s}$. We have obtained that $x$ belongs to the linear span of the vectors $e_1, \ldots, e_{k_{j_1}}$, and, as is easy to see, any vector from this space is orthogonal to the mixed system. Thus, the defect in this case is equal to $k_{j_1}$. In addition, each $k_j$ can be obtained as a defect. For this, it suffices to consider $\sigma=\{(s+1)n+j+1\}_{n=1}^\infty$. This completes the proof in the case of a finite set of integer defects. \medskip

\noindent  \textbf{Step 3: set of defects including infinity} \medskip

Similarly, we can construct an example of a system with a set of defects $\{0=k_0, k_1, k_2, \ldots, k_s\, \ldots\ , \infty \}$. We include here the case when the set of defects is given by $ \{k_0, k_1, ..., k_s, \infty \}$. In this case we assume that $k_j = k_s$ for $j\ge s+1$. Let $H_0$ be an infinite-dimensional closed subspace of $H$ that has an infinite-dimensional orthogonal complement. We denote by $\{f_j\}$ the orthonormal basis of $H_0^{\bot}$, and by $\{e_n\}$ the orthonormal basis of $H_0$. Put 
 $$
 x_n^s:= 2^n \sum \limits_{j=k_s +1}^n \frac{1}{n^{j-1}} f_j +e_n
 $$
 (see the construction in the Lemma \ref{example_inf}). Now we define a sequence $x_n$ on segments of the natural numbers from $n_m+1:=\frac{m(m+1)}{2} +1$ to $n_{m+1}=\frac{(m+1)(m+2)}{2} $ ($m=0,1,2, \ldots$). Let us define
 $$
 x_{n_m +1}:=x_{n_m +1}^{0}, \ x_{n_m +2}:=x_{n_m +2}^{1}, \ldots, \ x_{n_{m+1}}:=x_{n_{m+1}}^{m}.
 $$
 The sequence $x_n$ obeys the following simple pattern 
 $$
 x_1^{0}, \  x_2^{0}, \ x_3^{1}, \ x_4^{0}, \ x_5^{1}, \ x_6^{2}, \ldots
 $$
 It is clear that for any $m=0,1,2, \ldots$ among $x_n$ there are infinitely many vectors from the sequence $\{x^m_k\}$. Direct calculation shows that $\{e_n\}_{n=1}^\infty$ is biorthogonal to the constructed system. We will establish the completeness of the system $\{x_n\}_{n=1}^\infty$ by finding all possible defects of mixed systems. Let $\sigma$ be a subset of natural numbers and $\sigma=\sigma_{j_1}\sqcup \sigma_{j_2}\sqcup \ldots$, where $\sigma_{j_l}$ are sets such that $x_n=x_n^{j_l}$ on them. Let us find the defect of $\{x_n\}_{n=1}^\infty$ with respect to $\sigma$. First, let there be an infinite set among $\sigma_{j_l}$. Let $j_1$ be the smallest of those numbers $j$ for which $\sigma_j$ is infinite. Without loss of generality, among the elements of $\sigma$ there are no $n$ such that $x_n=x_n^m$ for $m<j_1$ (for each $m<j_1$ the number of such $n \in \sigma$ is finite by the definition of $j_1$, which means that they can all be excluded from $\sigma$, without changing the defect of the mixed system). Consider the vector $x\in H$ orthogonal to the vectors $x_n=x_n^{j_1}, n \in \sigma_{j_1}$. Let $x=\sum \limits_{n} c_n e_n + \sum \limits_{j} d_j f_j$. The orthogonality conditions are written as 
 \begin{equation} \label{ort_inf_1}
    2^n \sum \limits_{j=k_{j_1} +1}^n \frac{1}{n^{j-1}} d_j +c_n=0, \quad n\in \sigma_{j_1}.  
 \end{equation}
From here, as in the Lemma \ref{example_inf}, it can be proven that $d_j=0$ for $j\geq k_{j_1}+1$, and hence $c_n=0,n\in \sigma_{j_1}$. But for $n \in \sigma \setminus \sigma_{j_1}$ equality $x_n=x_n^m$, $m>j_1$, holds, and for such $x_n^m$ the orthogonality condition with vector $x$ takes the form
\begin{equation} \label{ort_inf_2}
 2^n \sum \limits_{j=k_{m} +1}^n \frac{1}{n^{j-1}} d_j +c_n=0.   
\end{equation}
Moreover, $k_m \geq k_{j_1}$, and therefore for these $n$ the equalities $c_n=0$ are fulfilled. For $n \notin \sigma$ the equality $c_n=0$ is obvious, since $x_n^*=e_n$. As a result, we have obtained that if $x$ is orthogonal to the mixed system, then it can be represented as $x=\sum \limits_{j=1}^{k_{j_1}} d_j f_j$.

It is easy to see that, conversely, any $x$ of this form satisfies all the required equations. From this we obtain that the defect of $x_n$ with respect to $\sigma$ is equal to $k_{j_1}$, and any $k_j$ is realized as a defect, since for $\sigma$ we can take exactly those numbers $n$ for which $x_n=x_n^{j}$ (in particular, for $\sigma=\mathbb{N}$ the set $\sigma_{j_1}$ from the proof coincides with those $n$ for which $x_n=x_n^{0}$, and hence any $x$ orthogonal to the entire system is zero). 

It remains to show that the defect is infinite when all $\sigma_{j_l}$ are finite. We will search for a vector $x$ orthogonal to the mixed system in the form $x=f_j + x'$, where $x' \in Lin\{e_n\}_{n=1}^\infty$. It is clear that for any $j$ such $x$ can be found since only for a finite number of $n \in \sigma$ the expansion of $x_n$ contains a multiple of $f_j$. For those $n$ from $\sigma$ for which the expansion of $x_n$ contains a multiple of $f_j$ the equality $c_n=-\frac{2^n}{n^{j-1}}$ must be satisfied by \eqref{ort_inf_2}. For the rest, $c_n=0$ also due to \eqref{ort_inf_2}, or due to the fact that $x_n^*=e_n$. Conversely, by defining $c_n$ in this way, we obtain a vector from the orthogonal complement to the mixed system. The family of vectors $x$ of the indicated type forms an infinite-dimensional space, and therefore, in the case under consideration, the defect with respect to $\sigma$ is infinite. 
\end{proof}

\end{document}